\documentclass[11pt]{amsart}
\usepackage{setspace}
\usepackage{a4}
\usepackage{amssymb,amsmath,amsthm,latexsym}

\usepackage[left=20mm,top=0.6in,bottom=12mm]{geometry}

\usepackage{amsfonts}
\usepackage{amsfonts}
\usepackage{graphicx}
\usepackage{textcomp}
\usepackage{cite}
\usepackage{enumerate}
\usepackage[mathscr]{euscript}
\usepackage{mathtools}
\newtheorem{theorem}{Theorem}[section]

\newtheorem{definition}[theorem]{Definition}
\newtheorem{example}[theorem]{Example}

\newtheorem{problem}[theorem]{Problem}
\newtheorem{proposition}[theorem]{Proposition}
\newtheorem{remark}[theorem]{Remark}
\title{This is the title}
\usepackage{amssymb}
\usepackage{amssymb}
\usepackage{amssymb}
\usepackage{amssymb}
\usepackage{amsmath}
\usepackage{tikz}
\usepackage{hyperref}
\usepackage{enumerate}
\usepackage{mathtools}
\usepackage{amsmath}
\usepackage{tikz}
\usepackage{amssymb}
\usepackage{amsmath}
\usepackage{tikz}
\usepackage{hyperref}
\begin{document}
	\title
	{{Approximately Dual \lowercase{p}-Approximate Schauder Frames}}
	\author{K. Mahesh Krishna}
	\address{$^1$Stat-Math Unit,	Indian Statistical Institute, Bangalore Centre,
		Karnataka 560 059 India}
	\email{kmaheshak@gmail.com}

	\author{P. Sam Johnson}
	\address{$^2$Department of Mathematical and Computational Sciences, National Institute of Technology Karnataka (NITK), Surathkal, Mangaluru 575 025, India}
	\email{sam@nitk.edu.in}
	
	\maketitle

\begin{abstract}  Difficulty in the   construction of  dual frames for a given  Hilbert space   led   to the introduction of approximately dual frames in Hilbert spaces by Christensen and Laugesen. It becomes even more difficult in Banach spaces to construct duals.  For this purpose, we introduce approximately dual frames for a class of approximate Schauder frames for Banach spaces and develop basic theory. Approximate duals for this subclass is completely characterized and its perturbation is also studied.
\end{abstract}




\section{Introduction}
Works of Holub \cite{HOLUB} and Li \cite{LI} classify frames and its duals for Hilbert spaces using the standard orthonormal basis for the standard separable Hilbert space. In the course of classifying approximate Schauder frames for Banach spaces using standard Schauder bases for classical sequence spaces, we \cite{KRISHNAJOHNSON} introduced the notion of p-approximate Schauder frames (p-ASFs) for Banach spaces (for approximate Schauder frames see \cite{FREEMANODELL, THOMAS}) which behaves in much more similar way that of frames for Hilbert spaces whose definition reads as follows. Let $\mathcal{X}$ be a separable Banach space and $\mathcal{X}^*$ be its dual. 
 \begin{definition}\cite{KRISHNAJOHNSON}\label{PASFDEF}
	An ASF $ (\{f_n \}_{n}, \{\tau_n \}_{n}) $  for $\mathcal{X}$	is said to be a \textbf{p-approximate Schauder frame} (p-ASF), $p \in [1, \infty)$ if the following conditions hold.
	\begin{enumerate}[\upshape(i)]
		\item The \textbf{frame operator} $	S_{f, \tau}:\mathcal{X}\ni x \mapsto S_{f, \tau}x\coloneqq \sum_{n=1}^\infty
		f_n(x)\tau_n \in
		\mathcal{X}$ is a well-defined bounded linear invertible eoperator. 
		\item The \textbf{analysis operator} $\theta_f: \mathcal{X}\ni x \mapsto \theta_f x\coloneqq \{f_n(x)\}_n \in \ell^p(\mathbb{N})$  is a well-defined bounded linear operator. 
		\item The \textbf{synthesis operator} $\theta_\tau : \ell^p(\mathbb{N}) \ni \{a_n\}_n \mapsto \theta_\tau \{a_n\}_n\coloneqq \sum_{n=1}^\infty a_n\tau_n \in \mathcal{X} $ is a well-defined bounded linear operator. 
	\end{enumerate}
Constants  $a,b>0$ satisfying 
\begin{align*}
a\|x\|\leq \left\|\sum_{n=1}^\infty
f_n(x)\tau_n \right\|\leq b\|x\|, \quad \forall x \in  \mathcal{X}
\end{align*} 
are called  as lower ASF bound and   upper ASF bound, respectively.
\end{definition}
Since the frame operator $S_{f, \tau}$ is invertible, we have 
\begin{align}\label{DUALMOT}
x=\sum_{n=1}^\infty(f_nS_{f, \tau}^{-1})(x)\tau_n=\sum_{n=1}^\infty f_n(x)S_{f, \tau}^{-1}\tau_n, \quad\forall x \in  \mathcal{X}.
\end{align}
It is also easy to see that both $ (\{f_nS_{f, \tau}^{-1}\}_{n}, \{\tau_n\}_{n}) $ and $ (\{f_n\}_{n}, \{S_{f, \tau}^{-1}\tau_n\}_{n}) $ are p-ASFs for $ \mathcal{X}$. In general, there are other p-ASFs that satisfy Equation (\ref{DUALMOT}). This leads to the following notion.
\begin{definition}\cite{KRISHNAJOHNSON}\label{DUALITYMINE}
	Let $ (\{f_n\}_{n}, \{\tau_n\}_{n}) $ be a p-ASF for 	$\mathcal{X}$. 	A p-ASF $ (\{g_n \}_{n}, \{\omega_n \}_{n}) $ for $\mathcal{X}$ is a \textbf{dual}  for $ (\{f_n \}_{n}, \{\tau_n \}_{n}) $ if 
	\begin{align*}
		x=\sum_{n=1}^\infty g_n(x) \tau_n=\sum_{n=1}^\infty
		f_n(x) \omega_n, \quad \forall x \in
		\mathcal{X}.
	\end{align*}
\end{definition}
 Following theorem classifies duals of  p-ASFs, where $\{e_n\}_n$ denotes the standard Schauder basis for  $\ell^p(\mathbb{N})$, $p \in [1, \infty)$ and $\{\zeta_n\}_n$ denotes the coordinate functionals associated with $\{e_n\}_n$.
 \begin{theorem}\label{ALLDUAL}\cite{KRISHNAJOHNSON}\
 Let $ (\{f_n \}_{n}, \{\tau_n \}_{n}) $ be a  p-ASF for   $\mathcal{X}$. Then a  p-ASF  $ (\{g_n \}_{n}, \{\omega_n \}_{n}) $ for $\mathcal{X}$ is a dual  for $ (\{f_n \}_{n}, \{\tau_n \}_{n}) $ if and only if
 \begin{align*}
 	&g_n=f_nS_{f,\tau}^{-1}+\zeta_nU-f_nS_{f,\tau}^{-1}\theta_\tau U,\\
 	&\omega_n=S_{f,\tau}^{-1}\tau_n+Ve_n-V\theta_fS_{f,\tau}^{-1}\tau_n, \quad \forall n \in \mathbb{N}
 \end{align*}
 such that the operator 
 \begin{align*}
 	S_{f,\tau}^{-1}+VU-V\theta_fS_{f,\tau}^{-1}\theta_\tau U
 \end{align*}
 is bounded invertible, where   $U:\mathcal{X} \to \ell^p(\mathbb{N})$ and $ V:\ell^p(\mathbb{N})\to \mathcal{X}$ are bounded linear operators.
\end{theorem}

It is known even in the case of Hilbert spaces that (see \cite{CHRISTENSENLAUGESEN}) it is difficult to construct or verify whether a given collection is a dual for a frame. This leads to the notion of approximately dual frames by Christensen and Laugensen \cite{CHRISTENSENLAUGESEN} which is motivated from the work of Li and Yan \cite{LIYAN}. Following the paper \cite{CHRISTENSENLAUGESEN}, there is a plenty of activity on the notion of approximately dual frames and its applications to wavelet frames, Gabor systems, shift-invarinat systems, localized frames, cross Gram matrices   \cite{DORFLERMATUSIAK, FEICHTINGERONCHISWIESMEYR, FEICHTINGERGRYBOSONCHIS, BENAVENTECHRISTENSENZAKOWICZ, LIANYOU, CHRISTENSENJANSSENKIMKIM, CHRISTENSENKIMKIM, BALAZSSHAMSABADIAREFIJAMAALRAHIMI} etc. This notion is also useful in the study of famous  Mexican hat problem \cite{BUILAUGESENMEXICAN1, BUILAUGESENMEXICAN2}. One can see Chapter 6 in the monograph \cite{CHRISTENSENBOOK} which gives a snapshot of approximately dual frames for Hilbert spaces. 
  
  The purpose of this paper is to study the notion of approximately dual frames for p-approximately Schauder frames for Banach spaces. One can naturally ask whether the notion of approximately dual frames can be defined for approximately Schauder frames. It seems that one can formulate the notion but it is difficult to study further due to the failure of factorization property of approximaly Schauder frames.

  \section{Approximately dual p-approximate Schauder frames}
 To define the notion of approximately dual frames for Hilbert spaces one need not consider frames, the notion of Bessel sequences suffices. To do the same thing for Banach spaces we first weaken Definition \ref{PASFDEF}. Our motivation comes from a characterization of Bessel sequences   for Hilbert spaces which reads as follows. 
  \begin{theorem}\cite{CHRISTENSENBOOK}\label{BESSELCHARACTERIZATION}
  	For a collection $\{\tau_n\}_n$ in a separable Hilbert space  $\mathcal{H}$, the following are equivalent. 
  		\begin{enumerate}[\upshape (i)]
  		\item $\{\tau_n\}_n$ is a Bessel sequence for   $\mathcal{H}$, i.e., there exist   $b>0$ such that 
$	\sum_{n=1}^\infty |\langle h, \tau_n\rangle|^2\leq b\|h\|^2,  \forall h \in \mathcal{H}.$
\item The map $ \mathcal{H}\ni h \mapsto  \{\langle h, \tau_n \rangle \}_n \in \ell^2(\mathbb{N})$  is a well-defined bounded linear operator. 
\item The map $ \ell^2(\mathbb{N}) \ni \{a_n\}_n \mapsto  \sum_{n=1}^\infty a_n\tau_n \in \mathcal{H}$ is a well-defined bounded linear operator.
  	\end{enumerate}
  	\end{theorem}
  Theorem \ref{BESSELCHARACTERIZATION}    leads to the following in Banach spaces.
  \begin{definition}
  	A pair  $ (\{f_n \}_{n}, \{\tau_n \}_{n}) $ is said to be a \textbf{p-approximate Bessel sequence}  for $\mathcal{X}$ (p-ABS), $p \in [1, \infty)$ if the following conditions hold.
  	\begin{enumerate}[\upshape(i)]
  	\item The \textbf{analysis operator} $\theta_f: \mathcal{X}\ni x \mapsto \theta_f x\coloneqq \{f_n(x)\}_n \in \ell^p(\mathbb{N})$  is a well-defined bounded linear operator. 
  		\item The \textbf{synthesis operator} $\theta_\tau : \ell^p(\mathbb{N}) \ni \{a_n\}_n \mapsto \theta_\tau \{a_n\}_n\coloneqq \sum_{n=1}^\infty a_n\tau_n \in \mathcal{X} $ is a well-defined bounded linear operator. 
  	\end{enumerate}
  Constants $c,d>0$ satisfynig 
  \begin{align*}
  \left(\sum_{n=1}^\infty
  |f_n(x)|^p\right)^\frac{1}{p}\leq c \|x\|, \quad \forall x \in \mathcal{X}, \quad 
  \left\|\sum_{n=1}^\infty a_n\tau_n\right\|\leq d \left(\sum_{n=1}^\infty
  |a_n|^p\right)^\frac{1}{p}, \quad \forall \{a_n\}_n  \in \ell^p(\mathbb{N}),
  \end{align*}
  are called as analysis and synthesis bounds, respectively.
  \end{definition}
  Following result gives examples of p-ABS and it also characterizes them.
   
   \begin{theorem}\label{THAFSCHAR}
  	A pair  $ (\{f_n\}_{n}, \{\tau_n\}_{n}) $ is a p-ABS for 	$\mathcal{X}$
  	if and only if 
  	\begin{align*}
  		f_n=\zeta_n U, \quad \tau_n=Ve_n, \quad \forall n \in \mathbb{N},
  	\end{align*}  
  	where $U:\mathcal{X} \rightarrow\ell^p(\mathbb{N})$, $ V: \ell^p(\mathbb{N})\to \mathcal{X}$ are bounded linear operators. 
  \end{theorem}
  \begin{proof}
  	$(\Leftarrow)$ It is easy to see that $\theta_f$ and $\theta_\tau$ are bounded linear operators. \\
  	$(\Rightarrow)$ Define $U\coloneqq \theta_f$, $V\coloneqq \theta_\tau$. Then $\zeta_nUx=\zeta_n\theta_fx=\zeta_n(\{f_k(x)\}_k)=f_n(x)$, $\forall x \in \mathcal{X}$, $Ve_n=\theta_\tau e_n=\tau_n$, $\forall n \in \mathbb{N}$. 
  \end{proof}
The idea behind the notion of approximate dual frames is to weaken the duality condition in Definition \ref{DUALITYMINE}.
 \begin{definition}\label{APPROXIMATELYDUALPABS}
 	Let  $ (\{f_n\}_{n}, \{\tau_n\}_{n}) $ and  $ (\{g_n\}_{n}, \{\omega_n\}_{n}) $ be two  p-ABSs for 	$\mathcal{X}$. We say that they are \textbf{approximately dual p-ABSs} if 
 	\begin{align*}
 	\|I_\mathcal{X}-\theta_\omega\theta_f\|< 1 \quad \text{ and } \quad \|I_\mathcal{X}-\theta_\tau\theta_g\|< 1.
 	\end{align*} 
 	In addition, if $ (\{f_n\}_{n}, \{\tau_n\}_{n}) $ and  $ (\{g_n\}_{n}, \{\omega_n\}_{n}) $ are   p-ASFs for 	$\mathcal{X}$, then we say that they are \textbf{approximately dual p-ASFs}.
 \end{definition} 
In the next example we see that both conditions in Definition  \ref{APPROXIMATELYDUALPABS}  are independent, this is contrast to the situation of Hilbert spaces in which adjoint operation reveals that one condition is enough.
\begin{example}\label{INDEPENDENT}
	Let $R$ be the right-shift and $L$ be the left-shift operator on $\ell^p(\mathbb{N})$. Define $f_n\coloneqq\zeta_n$, $\tau_n\coloneqq Re_n$, $g_n\coloneqq\zeta_nL$, $\omega_n\coloneqq e_n$. Then $ (\{f_n\}_{n}, \{\tau_n\}_{n}) $ and  $ (\{g_n\}_{n}, \{\omega_n\}_{n}) $ are   p-ABSs for $\ell^p(\mathbb{N})$ but 
	\begin{align*}
	\|I_{\ell^p(\mathbb{N})}-\theta_\omega\theta_f\|=\|I_{\ell^p(\mathbb{N})}-I_{\ell^p(\mathbb{N})}\|< 1, \quad 
	 \|I_{\ell^p(\mathbb{N})}-\theta_\tau\theta_g\|=\|I_{\ell^p(\mathbb{N})}-RL\|=1.
	\end{align*}
\end{example}
In the next example we show that there are p-ABSs which are not p-ASFs. This is in contrast with Hilbert space situation. Using Lemma 6.3.2 in \cite{CHRISTENSENBOOK} we can show that whenever two Bessel sequences for a Hilbert space are approximately duals, then they are necessarily frames.
\begin{example}\label{SECONDEXAMPLE}
	Let $R$  and $L$ be as in Example \ref{INDEPENDENT}. Define $f_n\coloneqq\zeta_n$, $\tau_n\coloneqq Le_n$, $g_n\coloneqq\zeta_nR$, $\omega_n\coloneqq e_n$. Then $ (\{f_n\}_{n}, \{\tau_n\}_{n}) $ and  $ (\{g_n\}_{n}, \{\omega_n\}_{n}) $ are   p-ABSs for $\ell^p(\mathbb{N})$ and 
	\begin{align*}
		\|I_{\ell^p(\mathbb{N})}-\theta_\omega\theta_f\|=\|I_{\ell^p(\mathbb{N})}-I_{\ell^p(\mathbb{N})}\|< 1, \quad 
		\|I_{\ell^p(\mathbb{N})}-\theta_\tau\theta_g\|=\|I_{\ell^p(\mathbb{N})}-LR\|<1.
	\end{align*}
Hence $ (\{f_n\}_{n}, \{\tau_n\}_{n}) $ and  $ (\{g_n\}_{n}, \{\omega_n\}_{n}) $ are   approximately dual p-ABSs for $\ell^p(\mathbb{N})$. Now note that $\theta_\tau\theta_f=L$ and  $\theta_\omega\theta_g=R$. Since both left and right shift operators are not invertible, both $ (\{f_n\}_{n}, \{\tau_n\}_{n}) $ and  $ (\{g_n\}_{n}, \{\omega_n\}_{n}) $ are   not  p-ASFs for $\ell^p(\mathbb{N})$.
\end{example}
\begin{remark}
	\begin{enumerate}[\upshape(i)]
		\item Inequalities in Definition \ref{APPROXIMATELYDUALPABS}  say that both $S_{f, \omega}$ and $S_{g, \tau}$ are invertible. Thus if $ (\{f_n\}_{n}, $ $ \{\tau_n\}_{n}) $ and  $ (\{g_n\}_{n}, \{\omega_n\}_{n}) $	are approximately dual p-ABSs for $\mathcal{X}$, then  $ (\{f_n\}_{n}, \{\omega_n\}_{n}) $ and  $ (\{g_n\}_{n}, $ $\{\tau_n\}_{n}) $	are  p-ASFs for $\mathcal{X}$.
		\item Proposition 2.13 in \cite{KRISHNAJOHNSON}  says that dual p-ASFs are always approximately dual p-ASFs. Since dual p-ASF always exists for every p-ASF, it follows that every p-ASF has approximate dual p-ASF.
	\end{enumerate}
\end{remark}
 Definition  \ref{APPROXIMATELYDUALPABS}  gives the  following proposition. 
  \begin{proposition}
  If $ (\{f_n\}_{n}, \{\tau_n\}_{n}) $ and  $ (\{g_n\}_{n}, \{\omega_n\}_{n}) $	are approximately dual p-ABSs for $\mathcal{X}$, then 
  \begin{align*}
  \left\|x-\sum_{n=1}^\infty
  f_n(x)\omega_n \right\|< \|x\| \quad \text{ and } \quad \left\|x-\sum_{n=1}^\infty
  g_n(x)\tau_n \right\|< \|x\|, \quad \forall x \in \mathcal{X}\setminus\{0\}.
  \end{align*}
  \end{proposition}
  \begin{proof}
  	Given $ x \in \mathcal{X}\setminus\{0\}$, 
  	\begin{align*}
 & \left\|x-\sum_{n=1}^\infty
  	f_n(x)\omega_n \right\|=\|I_\mathcal{X}x-\theta_\omega\theta_fx\|\leq \|I_\mathcal{X}x-\theta_\omega\theta_f\|\|x\|< \|x\|, \\
  &	  \left\|x-\sum_{n=1}^\infty
  	g_n(x)\tau_n \right\|=\|I_\mathcal{X}x-\theta_\tau\theta_gx\|\leq \|I_\mathcal{X}-\theta_\tau\theta_g\|\|x\|< \|x\|.
  	\end{align*}
  \end{proof}
Our next agenda  gives a complete description of p-ABSs and p-AFS's. For this purpose, we weaken Definition \ref{DUALITYMINE}.
\begin{definition}
	Let $ (\{f_n\}_{n}, \{\tau_n\}_{n}) $ be a p-ABS for 	$\mathcal{X}$. 	A p-ABS  $ (\{g_n \}_{n}, \{\omega_n \}_{n}) $ for $\mathcal{X}$ is a \textbf{dual}  for $ (\{f_n \}_{n}, \{\tau_n \}_{n}) $ if 
	\begin{align*}
		x=\sum_{n=1}^\infty g_n(x) \tau_n=\sum_{n=1}^\infty
		f_n(x) \omega_n, \quad \forall x \in
		\mathcal{X}.
	\end{align*}
\end{definition}
In the case of Hilbert spaces, existence of dual Bessel sequence forces the Bessel sequence to a frame (Lemma 6.3.2 in \cite{CHRISTENSENBOOK}) However, Example \ref{SECONDEXAMPLE} reveals that existence of dual Bessel sequences do not imply that they are frames.
In the next proposition we show that approximate duals generate duals whose proof is a routine computation.
\begin{proposition}
Let  $ (\{f_n\}_{n}, \{\tau_n\}_{n}) $ and  $ (\{g_n\}_{n}, \{\omega_n\}_{n})$ be approximately dual   p-ABSs for 	$\mathcal{X}$. Then $ (\{g_nS_{g,\tau}^{-1}\}_{n}, \{S_{f,\omega}^{-1}\omega_n\}_{n}) $ is a dual p-ABS for $ (\{f_n\}_{n}, \{\tau_n\}_{n}) $ 	and $ (\{f_nS_{f,\omega}^{-1}\}_{n}, \{S_{g,\tau}^{-1}\tau_n\}_{n}) $ is a dual p-ABS for $ (\{g_n\}_{n}, \{\omega_n\}_{n})$. 
\end{proposition}
An operator-theoretic description of approximately duals can be given which is illustrated in the next result.
\begin{theorem}\label{APPROXIMATEBESSELCHAR}
	Let $ (\{f_n\}_{n}, \{\tau_n\}_{n}) $ be a p-ABS for 	$\mathcal{X}$. 	A p-ABS  $ (\{g_n \}_{n}, \{\omega_n \}_{n}) $ for $\mathcal{X}$ is an approximately dual p-ABS   for $ (\{f_n \}_{n}, \{\tau_n \}_{n}) $ if and only if there exist bounded linear operators  $U,V:	\mathcal{X} \to \mathcal{X}$ satisfying 	$\|I_\mathcal{X}-U\|< 1$ and $\|I_\mathcal{X}-V\|< 1$ such that $ (\{h_n\coloneqq g_nU^{-1}\}_{n}, \{\rho_n\coloneqq V^{-1}\omega_n\}_{n}) $ is a dual for $ (\{f_n\}_{n}, \{\tau_n\}_{n}) $. Statement holds even if p-ABS is replaced by p-ASF.
\end{theorem}
\begin{proof}
$(\Rightarrow)$ Define $U\coloneqq \theta_\tau\theta_g$ and $V\coloneqq \theta_\omega\theta_f$. We then have $\|I_\mathcal{X}-U\|< 1$,  $\|I_\mathcal{X}-V\|< 1$ and 
\begin{align*}
	&\sum_{n=1}^\infty h_n(x) \tau_n=\sum_{n=1}^\infty g_n(U^{-1}x) \tau_n=\theta_\tau\theta_gU^{-1}x=UU^{-1}x=x, \\
	&\sum_{n=1}^\infty f_n(x) \rho_n=\sum_{n=1}^\infty
	f_n(x) V^{-1}\omega_n=V^{-1}\theta_\omega\theta_f=V^{-1}Vx=x, \quad \forall x \in
	\mathcal{X}.
\end{align*}
Hence $ (\{h_n\}_{n}, \{\rho_n\}_{n}) $ is a dual for $ (\{f_n\}_{n}, \{\tau_n\}_{n}) $.

$(\Leftarrow)$ A direct computation says that $\theta_h=\theta_gU^{-1}$ and $\theta_\rho=V^{-1}\theta_\omega$. Since $ (\{h_n\}_{n}, \{\rho_n\}_{n}) $ is a dual for $ (\{f_n\}_{n}, \{\tau_n\}_{n}) $, we get that 
\begin{align*}
&\|I_\mathcal{X}-\theta_\omega\theta_f\|	=\|I_\mathcal{X}-V\theta_\rho\theta_f\|=\|I_\mathcal{X}-V\|<1,\\
&\|I_\mathcal{X}-\theta_\tau\theta_g\|	=\|I_\mathcal{X}-\theta_\tau\theta_hU\|=\|I_\mathcal{X}-U\|<1.
\end{align*}
Hence $ (\{g_n \}_{n}, \{\omega_n \}_{n}) $  is an approximately dual p-ABS   for $ (\{f_n \}_{n}, \{\tau_n \}_{n}) $.
\end{proof}
Theorem  \ref{APPROXIMATEBESSELCHAR}  when combined with Theorem  \ref{ALLDUAL}   gives the following result for approximate dual frames.
\begin{theorem}
	Let $ (\{f_n\}_{n}, \{\tau_n\}_{n}) $ be a p-ASF for 	$\mathcal{X}$. 	A p-ASF  $ (\{g_n \}_{n}, \{\omega_n \}_{n}) $ for $\mathcal{X}$ is an approximately dual p-ASF   for $ (\{f_n \}_{n}, \{\tau_n \}_{n}) $ if and only if there exist bounded linear operators  $U,V:	\mathcal{X} \to \mathcal{X}$, $A:\mathcal{X} \to \ell^p(\mathbb{N})$, $ B:\ell^p(\mathbb{N})\to \mathcal{X}$ satisfying 	$\|I_\mathcal{X}-U\|< 1$, $\|I_\mathcal{X}-V\|< 1$, 
	\begin{align*}
		&g_n=f_nS_{f,\tau}^{-1}U+\zeta_nAU-f_nS_{f,\tau}^{-1}\theta_\tau AU,\\
		&\omega_n=VS_{f,\tau}^{-1}\tau_n+VBe_n-VB\theta_fS_{f,\tau}^{-1}\tau_n, \quad \forall n \in \mathbb{N}
	\end{align*}
	such that the operator 
	\begin{align*}
		S_{f,\tau}^{-1}+BA-V\theta_fS_{f,\tau}^{-1}\theta_\tau A
	\end{align*}
	is bounded invertible.
	\end{theorem}
In \cite{CHRISTENSENLAUGESEN},  Christensen and Laugensen constructed  iterations which construct approximately duals at each stage and converge to identity operator in operator norm. Here we have a  similar result for p-ABSs.
\begin{theorem}
Let 	$ (\{f_n\}_{n}, \{\tau_n\}_{n}) $ be a p-ABS for 	$\mathcal{X}$ and $ (\{g_n\}_{n}, \{\omega_n\}_{n}) $ be a p-ABS for 	$\mathcal{X}$ which is approximately dual for $ (\{f_n\}_{n}, \{\tau_n\}_{n}) $.
\begin{enumerate}[\upshape(i)]
	\item The dual p-ABS   $ (\{g_nS_{g,\tau}^{-1}\}_{n}, \{S_{f,\omega}^{-1}\omega_n\}_{n}) $  for $ (\{f_n\}_{n}, \{\tau_n\}_{n}) $ can be written as 
	\begin{align*}
	&g_nS_{g,\tau}^{-1}	=g_n+\sum_{m=1}^{\infty}g_n(I_\mathcal{X}-S_{g,\tau})^m,\\
	&S_{f,\omega}^{-1}\omega_n=\omega_n+\sum_{m=1}^{\infty}(I_\mathcal{X}-S_{f, \omega})^m\omega_n, \quad \forall n \in \mathbb{N}.
	\end{align*}
	\item For $N\in \mathbb{N}$, define 
	\begin{align*}
		&h_n^{(N)}\coloneqq  g_n+\sum_{m=1}^{N}g_n(I_\mathcal{X}-S_{g,\tau})^m,\\
		&\rho_n^{(N)}\coloneqq \omega_n+\sum_{m=1}^{N}(I_\mathcal{X}-S_{f, \omega})^m\omega_n, \quad \forall n \in \mathbb{N}.
	\end{align*}
Then $ (\{h_n^{(N)}\}_{n}, \{\tau_n^{(N)}\}_{n}) $ is a p-ABS and is dual for $ (\{f_n\}_{n}, \{\tau_n\}_{n}) $, for each $N\in \mathbb{N}$. Moreover, 
\begin{align*}
	&\|I_\mathcal{X}-\theta_\rho^{(N)}\theta_f\|\leq \|I_\mathcal{X}-S_{f,\omega}\|^{N+1}\to 0 \quad \text{ as } N \to \infty, \\
	&	\|I_\mathcal{X}-\theta_\tau\theta_h^{(N)}\|\leq \|I_\mathcal{X}-S_{g,\tau}\|^{N+1}\to 0 \quad \text{ as } N \to \infty.
\end{align*}
\end{enumerate}
\end{theorem}
\begin{proof}
	\begin{enumerate}[\upshape(i)]
		\item This follows from the Neumann series 
		\begin{align*}
			S_{g,\tau}^{-1}=\sum_{m=0}^{\infty}(I_\mathcal{X}-S_{g,\tau})^m, \quad S_{f,\omega}^{-1}=\sum_{m=0}^{\infty}(I_\mathcal{X}-S_{f,\omega})^m.
		\end{align*}
		\item Clearly $ (\{h_n^{(N)}\}_{n}, \{\tau_n^{(N)}\}_{n}) $ is a p-ABS. Consider 
		
		Hence $ (\{h_n^{(N)}\}_{n}, \{\tau_n^{(N)}\}_{n}) $ is a dual  for $ (\{f_n\}_{n}, \{\tau_n\}_{n}) $. Now consider
		\begin{align*}
		\theta_\rho^{(N)}\theta_fx&=\sum_{n=1}^{\infty}f_n(x)\rho_n^{(N)}=\sum_{n=1}^{\infty}f_n(x)\sum_{m=0}^{N}(I_\mathcal{X}-S_{f, \omega})^m\omega_n\\
		&=\sum_{m=0}^{N}(I_\mathcal{X}-S_{f, \omega})^m\sum_{n=1}^{\infty}f_n(x)\omega_n=\sum_{m=0}^{N}(I_\mathcal{X}-S_{f, \omega})^mS_{f,\omega}x\\
		&=\sum_{m=0}^{N}(I_\mathcal{X}-S_{f, \omega})^m(I_\mathcal{X}-(I_\mathcal{X}-S_{f,\omega}))x=x-(I_\mathcal{X}-S_{f, \omega})^{N+1}x,\quad \forall x \in \mathcal{X}\\
		&\implies \|I_\mathcal{X}-\theta_\rho^{(N)}\theta_f\|\leq \|I_\mathcal{X}-S_{f,\omega}\|^{N+1}
		\end{align*}
	and 
	\begin{align*}
	\theta_\tau\theta_h^{(N)}x&=	\sum_{n=1}^{\infty}h_n^{(N)}(x)\tau_n=\sum_{n=1}^{\infty}\sum_{m=0}^Ng_n((I_\mathcal{X}-S_{g,\tau})^mx)\tau_n\\
	&=\sum_{m=0}^N\sum_{n=1}^\infty g_n((I_\mathcal{X}-S_{g,\tau})^mx)\tau_n=\sum_{m=0}^NS_{g,\tau}(I_\mathcal{X}-S_{g,\tau})^mx\\
	&=\sum_{m=0}^N(I_\mathcal{X}-(I_\mathcal{X}-S_{g,\tau}))(I_\mathcal{X}-S_{g,\tau})^mx=x-(I_\mathcal{X}-S_{g,\tau})^{N+1}x, \quad  \forall x \in \mathcal{X}\\
	&\implies \|I_\mathcal{X}-\theta_\tau\theta_h^{(N)}\|\leq \|I_\mathcal{X}-S_{g,\tau}\|^{N+1}.
	\end{align*}
Conclusion follows from $\|I_\mathcal{X}-S_{f,\omega}\|<1$ and $\|I_\mathcal{X}-S_{g,\tau}\|<1$.
	\end{enumerate}
\end{proof}
  In \cite{CHRISTENSENLAUGESEN}  Christensen and Laugensen showed that by considering sequences close to a given frame one can generate approximate duals. Here is the similar result in the context of Banach spaces.
  \begin{theorem}
  	Let $\{\tau_n\}_n$ be a collection in $\mathcal{X}$ and $\{f_n\}_n$ be a collection in $\mathcal{X}^*$.  	Let  $ (\{h_n\}_{n}, \{\rho_n\}_{n}) $ be a p-ASF for 	$\mathcal{X}$ such that 
  	\begin{align}\label{ANALYSISEST}
  	\left(\sum_{n=1}^{\infty}|f_n(x)-h_n(x)	|^p\right)^\frac{1}{p}\leq R\|x\|,\quad \forall x \in \mathcal{X}
  	\end{align}
  and 
  	\begin{align}\label{SYNTHESISEST}
 \left\|\sum_{n=1}^{\infty}a_n(\rho_n-\tau_n)\right\| 	\leq Q\left(\sum_{n=1}^{\infty}|a_n	|^p\right)^\frac{1}{p},\quad  \forall \{a_n\}_n  \in \ell^p(\mathbb{N}).
  \end{align}
  Let $ (\{g_n\}_{n}, \{\omega_n\}_{n}) $ be a dual frame for $ (\{h_n\}_{n}, \{\rho_n\}_{n}) $ with analysis bound $c$  and synthesis bound $d$. If $dR,cQ<1$, 
  then  $ (\{g_n\}_{n}, \{\omega_n\}_{n}) $ is an approximately dual  for $ (\{f_n\}_{n}, \{\tau_n\}_{n}) $.
  \end{theorem}
\begin{proof}
Inequalities (\ref{ANALYSISEST}) and (\ref{SYNTHESISEST})  (using Minkowski's and triangle inequalities) say that $ (\{f_n\}_{n}, \{\tau_n\}_{n}) $ is a p-ABS for $\mathcal{X}$.	Note that Inequality (\ref{ANALYSISEST}) can be written as $\|\theta_hx-\theta_fx\|\leq R\|x\|$, $\forall x \in \mathcal{X}$. Similarly Inequality (\ref{SYNTHESISEST}) can be written as $\|\theta_\tau\{a_n\}_n-\theta_\rho\{a_n\}_n\|\leq Q\|\{a_n\}_n\|, \forall \{a_n\}_n  \in \ell^p(\mathbb{N})$. Using these, we have 
\begin{align*}
	\|I_{\ell^p(\mathbb{N})}-\theta_\omega\theta_f\|=\|\theta_\omega\theta_h-\theta_\omega\theta_f\|\leq \|\theta_\omega\|\|\theta_h-\theta_f\|\leq dR< 1, 
\end{align*}	
and 
\begin{align*}
	 \|I_{\ell^p(\mathbb{N})}-\theta_\tau\theta_g\|=\|\theta_\rho\theta_g-\theta_\tau\theta_g\|\leq =\|\theta_\rho-\theta_\tau\|\|\theta_g\|\leq cQ<1
\end{align*}
which is required.
\end{proof}

\section{An Open Problem}
  To motivate the reader for further study on approximately duals of p-ASFs, we end the paper with an open problem. Our motivation comes from Hilbert space frame theory. Recall that  \cite{BALANCASAZZAHEIL} given a frame $ \{\tau_n\}_{n}$ for a Hilbert space $\mathcal{H}$, the \textbf{excess} of  $\{\tau_n\}_{n}$ is defined as 
  \begin{align*}
  		\text{Exc}(\{\tau_n\}_{n})\coloneqq \sup\{|\mathbb{M}|:\mathbb{M}\subseteq   \mathbb{N} ,\overline{\text{span}}\{\tau_n\}_{n\in \mathbb{N}\setminus\mathbb{M}}=\mathcal{H} \}.
  \end{align*}
  The following result shows that excess is an invariant for approximately dual frames.
  \begin{theorem}\cite{BAKICBERIC}
  	If $ \{\tau_n\}_{n}$ and $ \{\omega_n\}_{n}$ are approximately dual frames for a Hilbert space $\mathcal{H}$, then 
  	\begin{align*}
  		\text{Exc}(\{\tau_n\}_{n})=	\text{Exc}(\{\omega_n\}_{n}).	
  	\end{align*}
  \end{theorem}
We now formulate the notion of excess of p-ASFs and ask a problem.
\begin{definition}
Given a p-ASF	$ (\{f_n\}_{n}, \{\tau_n\}_{n}) $ for $\mathcal{X}$, we define the \textbf{p-excess} of $ (\{f_n\}_{n}, \{\tau_n\}_{n}) $  as 
\begin{align*}
\text{p-Exc}(\{f_n\}_{n}, \{\tau_n\}_{n})\coloneqq	\sup\{|\mathbb{M}|:\mathbb{M}\subseteq   \mathbb{N} ,\overline{\text{span}}\{\tau_n\}_{n\in \mathbb{N}\setminus\mathbb{M}}=\mathcal{X} , \overline{\text{span}}\{f_n\}_{n\in \mathbb{N}\setminus\mathbb{M}}=\mathcal{X}^*\}.
\end{align*}
\end{definition}
  \begin{problem}
  	Let $ (\{f_n\}_{n}, \{\tau_n\}_{n}) $ and  $ (\{g_n\}_{n}, \{\omega_n\}_{n}) $ be two  approximately dual p-ASFs for 	$\mathcal{X}$. Is the following relation
  	\begin{align*}
  		\text{p-Exc}(\{f_n\}_{n}, \{\tau_n\}_{n})=	\text{p-Exc}(\{g_n\}_{n}, \{\omega_n\}_{n})
  	\end{align*}
  	true?
  \end{problem}

 \bibliographystyle{plain}
 \bibliography{reference.bib}

\end{document}